\documentclass[letterpaper, 10pt, conference]{ieeeconf}

\IEEEoverridecommandlockouts                              
\usepackage{amsmath} 
\usepackage{amssymb}  
\usepackage{booktabs} 
\usepackage{graphicx}

\usepackage {amsfonts}

\newtheorem{lemma}{Lemma}

\newtheorem{remark}{Remark}

\newtheorem{definition}{Definition}
\newtheorem{problem}{Problem}

\newcommand{\N}{\mathbb{N}}
\newcommand{\R}{\mathbb{R}}
\def\è{\mbox{\`e}}
\def\ò{\mbox{\`o}}
\def\à{\mbox{\`a}}
\def\ì{\mbox{\`i}}
\def\ù{\mbox{\`u}}
\def\é{\mbox{\'e}}

\def\Pibf{\boldsymbol{\Pi}}

\def\urltilda{\kern -.15em\lower .7ex\hbox{\~{}}\kern .04em}

\title{\LARGE \bf
Two-step Nonnegative Matrix Factorization Algorithm for the Approximate Realization of Hidden Markov Models
}

\author{Lorenzo Finesso, Angela Grassi and Peter Spreij
\thanks{This work was partially supported by the PRIL project -- Regione Veneto}
\thanks{L. Finesso is with the Institute of Biomedical Engineering, ISIB-CNR, Padova, Italy.
        {\tt\footnotesize lorenzo.finesso@isib.cnr.it}}%
\thanks{A. Grassi is with the Department of Information Engineering, University of Padova, and a Research Associate with
        the Institute of Biomedical Engineering, ISIB-CNR, Padova, Italy.
        {\tt\footnotesize angela.grassi@dei.unipd.it}}%
\thanks{P. Spreij is with the University of Amsterdam, KdV Institute, The Netherlands.
        {\tt\footnotesize spreij@science.uva.nl}}%
}

\begin{document}

\maketitle
\thispagestyle{empty}
\pagestyle{empty}

%
%
%
%
\begin{abstract}

We propose a two-step algorithm for the construction of a Hidden Markov Model (HMM) of assigned size, i.e.
cardinality of the state space of the underlying Markov chain, whose $n$-dimensional
distribution is closest in divergence to a given distribution. The algorithm is
based on the factorization of a pseudo Hankel matrix, defined in terms of the given distribution,
into the product of a tall and a wide nonnegative matrix.
The implementation is based on the nonnegative matrix factorization (NMF) algorithm.
To evaluate the performance of our algorithm we produced some numerical simulations in the context
of HMM order reduction.

\end{abstract}


\section{Introduction}
Hidden Markov models (HMM) are a simple, yet very rich, class of
stochastic processes which has become ubiquitous in several areas of
signals, systems, and control. Here we restrict attention to HMMs
$(Y_t), t=1,2, \dots$ taking values into a finite set $\mathcal Y$.
The HMMs of this type can always be represented as deterministic functions of some
Markov chain $(S_t), t=1,2,\dots$, taking values in a finite set
$\mathcal S$, in the sense that $Y_t \stackrel{d}= f(S_t)$, for all
$t\in \N$, where the symbol $\stackrel{d}=$ means that the processes
$(Y_t)$ and $(f(S_t))$ have the same distributions. The cardinality
of $\mathcal S$ is called the size of the representation $Y_t
\stackrel{d}= f(S_t)$.

Since functions of Markov chains generally loose the Markov property,
HMMs can be used to model also dynamical behaviors exhibiting complex
dependencies from the past, which cannot be described within the class
of Markov chains. At the same time HMMs admit a simple parametric
description, through the transition probabilities matrix of $(S_t)$
and the deterministic function $f$. Being flexible and easy to describe
it comes as no surprise that the class
of HMMs is extensively employed in many applications to real data.
The applications  span the fields of engineering (modeling stochastic automata,
for automatic speech recognition and for communication networks), genetics (sequence
analysis), biology (to study neuro-transmission through
ion-channels), mathematical finance (to model rating transitions, or
to solve asset allocation problems), and many others.

Although by now many approaches to solve inferential problems about HMMs have been
proposed in the engineering and statistical literature, an algorithm to solve the exact stochastic
realization problem from distributional data is still missing.
Even the weakest form of stochastic realization, \emph{given the
finite dimensional distributions $p_Y(Y_1^n)$ of a HMM find the
parameters of any of its representations $(f(S_t))$}, has not yet been solved
satisfactorily.
The early attack dates back to~\cite{furst}, and~\cite{heller1965}, the
stochastic realization point of view was introduced later
in~\cite{picci1978} and~\cite{piccivans}, the most recent results
can be found in~\cite{bdo1999} and~\cite{vid2004}.

In the present paper we focus on the simpler \emph{approximate}
realization problem which can be roughly formulated as follows.
Given the distributions $q_Y(Y_1^n)$ of a stationary process, find a realization
of a HMM of assigned size, which best approximates $q_Y(Y_1^n)$ in divergence
rate, a most natural criterion of closeness between distributions.
Unfortunately there exists no general, closed form, analytic
expression of the divergence rate between
HMMs~\cite{karan},~\cite{hanmarcus2006}, let alone of the divergence rate between
a general stationary process and a HMM.
To obviate the difficulty
we formulate an alternative criterion, in terms of the informational
divergence between nonnegative, pseudo Hankel matrices, representing the finite
dimensional distributions of the processes. The approximate
realization problem becomes then amenable to the use of Nonnegative
Matrix Factorization (NMF) techniques. This approach was already
investigated in~\cite{finessospreij2002}, and~\cite{FGS-arxiv},
where we proposed a three step, NMF based, optimization procedure to construct
the parameters of the best approximate realization. The
same approach, in slightly different contexts, has been later
followed in~\cite{DLW} and~\cite{cybenko},
the latter proposing a two step algorithm based on estimated data
instead of exact distributional data, with the second step carried out via
linear programming.

The remainder of the paper is organized as follows. Section~\ref{sec:preliminaries}
contains preliminaries on HMMs. In Section~\ref{sec:hankel} we introduce the pseudo Hankel matrix of
the finite dimensional distributions of a stationary process
and study its factorization properties for the class of HMMs. Section~\ref{sec:algorithm}
introduces the two-step approximation algorithm. The last
Section discusses some of the numerical issues and provides examples
of order reduction for HMMs.


\section{Preliminaries}
\label{sec:preliminaries}
A stochastic process $(Y_{t})_{t\in\N}$ with values in ${\cal Y} =
\{1, 2,\dots m\}$ is a stationary \emph{hidden Markov model} (HMM)
of size $N$ if for some process $(X_{t})_{t\in \N}$, with values in
$\mathcal X = \{1, \dots N\}$, the pair $(X_t,Y_t)$ is jointly
stationary and, for all $y \in\mathcal{Y}$, $j \in \mathcal X$,
\begin{equation}\label{splitting property}
\begin{split} P(Y_{t+1}=y, X_{t+1}=j|X_t^-,Y_t^-) = \\ P(Y_{t+1}=y,
X_{t+1}=j|X_t).
\end{split}
\end{equation}
From~(\ref{splitting property}) it follows that both, the pair
$(S_t):=(X_t, Y_t)$, and $(X_t)$ alone, are Markov chains, a
property that shows why a HMM can always be described as a
deterministic function of a Markov chain, i.e. the projection onto
the second component of $(S_t)$. The distribution functions,
$p_Y(y_1\dots y_n)=P(Y_1^n =y_1^n)$, of $(Y_t)$ are specified by the
$m$ nonnegative $N \times N$ matrices
\begin{equation*}
m_{ij}(y)=P(Y_{t+1} = y, X_{t+1} = j \mid X_{t} = i),
~\label{m_ij(y)}
\end{equation*}
and by a row vector $\pi \! = \! \pi A$, where $A \! = \! \sum_{y}
M(y)$ is the transition probabilities matrix of~$(X_t)$. For any
(finite) string $w=y_1\cdots y_n$ one gets
$$p_Y(w) = \pi M(y_{1})
\cdots M(y_{n}) e,$$ where $e=(1,\ldots,1)^\top \in \R^N$. Define
$M(w)=M(y_1^n) = M(y_1)M(y_2)\dots M(y_n)$. It follows that, for any
pair of (finite length) strings $u$ and $v$, denoting by $uv$ their
concatenation, one has
\begin{equation}\label{p(uv)}
 p_Y(uv) =\pi M(u)  M(v) e =: \pi(u) \gamma(v).
\end{equation}
\indent For any finite length strings $u$ and $v$ and any $y \in \mathcal Y$
the vectors $\pi(u) \!:= \!\pi M(u)$ (row) and $\gamma(v) \!:=
\!M(v)e$ (column) satisfy the following relations
\begin{alignat}{2} \label{recs}
\pi(uy) & = \pi(u)M(y), & \qquad \pi(u) &=\sum_{y=1}^m \pi(yu), \nonumber \\
\gamma(yv) &= M(y)\gamma(v), & \qquad \gamma(v) &= \sum_{y=1}^m \gamma(vy)
\end{alignat}
In common usage the recurrence equations (the ones on the left) are called respectively
\emph{forward} equation (for $\pi$), and \emph{backward} equation (for $\gamma$).
%
%
%
%
\section{Hankel matrices and their factorizations}
\label{sec:hankel}
We introduce two lexicographic orders on the set $\mathcal{Y}^{n}$.
The \emph{first lexical order (flo)}, increasing from right to left,
and the \emph{last lexical order (llo)}, increasing from left to
right. By way of example, if $\mathcal{Y}=\{0,1\}$ and $n=2$, the
strings in increasing \emph{flo} are $(00,10,01,11)$, while, in
increasing \emph{llo}, are $(00,01,10,11)$.

\smallskip
\begin{definition} (Hankel matrices of stationary processes) If
$q_Y(\cdot)$ is the finite dimensional distribution of a stationary measure $Q$,
$$
\mathbf{H}^Q_{nn} := ||q_Y(u_r v_s)||_{r,s},
$$
where $u_r$ and $v_s$ run through \emph{all of} $\mathcal{Y}^{n}$ in \emph{flo} and in \emph{llo} respectively.
\end{definition}

\smallskip \noindent
The matrix $\mathbf{H}^Q_{nn}$ is generally called the Hankel matrix associated
to $Q$, see e.g.~\cite{bdo1999}, but it only has a structural resemblance to a
genuine Hankel matrix. The Hankel matrices of HMM measures admit several nonnegative
factorizations stemming from~(\ref{p(uv)}).
Specifically, if $p_Y$ is the finite dimensional distribution of a stationary HMM measure $P$,
\begin{eqnarray} \label{fatt blocco H_{KL}}
\mathbf{H}^P_{nn} &=&  \left[
\begin{array}{c} \pi(u_1) \nonumber \\
\vdots\\
\pi(u_{\ell})
\end{array}\right]
\left[
\begin{array}{ccc}\gamma(v_1)&\cdots
&\gamma(v_{\ell})
\end{array}\right] \\&=:& \boldsymbol{\Pi}_n\mathbf{\Gamma}_n,
\end{eqnarray}
where $\boldsymbol{\Pi}_n$, and $\mathbf{\Gamma}_n$
are matrices of sizes $m^n \times N$ and $N\times m^n$ respectively $(\ell = m^n).$

\noindent The matrices $\boldsymbol{\Pi}_n$, and $\mathbf{\Gamma}_n$, in turn, can
be factored as follows
\begin{equation} \label{Q1_L+1}
\mathbf{\Pi}_n= \left[\mathbf{\Pi}_{n-1}M(1);\ldots ;\mathbf{\Pi}_{n-1}M(m)\right],
\end{equation}
\begin{equation}
\mathbf{\Gamma}_{n}=\left[M(1)\mathbf{\Gamma}_{n-1}, \ldots ,M(m)\mathbf{\Gamma}_{n-1}\right],
\end{equation}
where '';'' and '','' denote, respectively, vertical and horizontal
stacking of the blocks. Introducing the matrix
$$
\mathbf{\widetilde{\Pi}}_{n} := \left[\mathbf{\Pi}_{n-1}M(1),
\ldots , \mathbf{\Pi}_{n-1}M(m)\right] \in \R_+^{m^{n-1}\times mN},
$$
equations~(\ref{Q1_L+1}) can be compactly written
\begin{equation} \label{factrecur}
\mathbf{\widetilde{\Pi}}_{n}=\mathbf{\Pi}_{n-1} \mathbf{M}, \qquad
\mathbf{\Gamma}_{n}=\mathbf{M}\mathbf{\Gamma}_{(n-1)},
\end{equation}
where
\begin{eqnarray}
\mathbf M  &:=& \left[M(y_1), \ldots , M(y_m)\right] \in \R_+^{N
\times mN}   \nonumber     \\
\mathbf{\Gamma}_{(n-1)}  &:=& {\rm diag}(\Gamma_{n-1} \dots
\Gamma_{n-1}) \in \R_+^{mN\times m^{n}}. \label{newgamma}
\end{eqnarray}

\smallskip
\begin{remark} \label{remark1} Note that both $\mathbf \Pi_{n-1}$ and $\mathbf \Gamma_{n-1}$
are readily obtained from $\mathbf \Pi_{n}$ and $\mathbf \Gamma_{n}$
by way of the marginal relations (the ones on the right) given in~(\ref{recs}).
\end{remark}


\section{Approximation problem and two-step algorithm}
\label{sec:algorithm}
We now pose the problem of best approximation of a given stationary
process with a HMM of assigned size, convert it into an approximate
nonnegative matrix factorization (NMF),\footnote{It has become commonplace in
the applied literature to reserve the
acronym NMF to denote the \emph{approximate} Nonnegative Matrix
Factorization problem, made popular by~\cite{lee}.} and propose an
algorithm for its numerical solution. The following definition will
be used throughout the paper.
\smallskip
\begin{definition} (Divergence between positive matrices) Let
$\mathbf{M},\mathbf{N}\in \mathbb{R}_+^{m\times n}$,
\begin{equation}
D(\mathbf{M}\|\mathbf{N}):=\sum_{ij}(M_{ij}\log\frac{M_{ij}}{N_{ij}}-M_{ij}+N_{ij}).~\label{def:I-div}
\end{equation}
\end{definition}
We are now ready to pose the approximation problem as a nonnegative
matrix factorization (NMF) problem.
\smallskip
\begin{problem} \label{approx-prob} Given a stationary
probability measure $Q$ on $\mathcal{Y}^\infty$ and $N \in
\mathbb{N}$, find the parameters $\{ M^*(y), \, y \in \mathcal{Y}
\}$ of a HMM of size $N$ whose Hankel matrix
$\mathbf{H}^*_{nn}=\boldsymbol\Pi_n^* \boldsymbol{\Gamma}_n^*$ ($n >
2N$) is closest to $\mathbf{H}^Q_{nn}$ in divergence. More
compactly: solve the minimization problem
\begin{equation} \label{fctpbm}
\underset{\mathbf \Pi_n,\mathbf{\Gamma}_n}\min D(\mathbf
H^Q_{nn}\|\boldsymbol\Pi_n \boldsymbol{\Gamma}_n),
\end{equation}
under constraints $\Pibf_n, \boldsymbol{\Gamma}_n \ge 0, \,\, e^\top
\Pibf_n e=1$, and $\boldsymbol{\Gamma}_n e=e$.
\end{problem}

\smallskip
The motivation for this setup comes from the fact that the
distributions of a HMM of size $N$ are fully determined by any of
its Hankel matrices $\mathbf{H}^P_{nn}$ with $n >2N$ see
e.g.~\cite{carlyle}. The constraints are imposed by the definitions
of the factors $\Pibf_n$ and $\boldsymbol{\Gamma}_n$ given above.

\smallskip A minimizing nonnegative factorization
$(\mathbf \Pi^*_n, \mathbf{\Gamma}^*_n)$ always exists,
see~\cite{finessospreij2005}, Proposition 2.1, but
Problem~\ref{approx-prob} also calls for the construction of the
corresponding parameters $M^*(y)$. The analysis of the ideal case
will serve as a guide. If $Q$ were a HMM law, the following
\emph{exact} nonnegative factorizations would hold by
equations~(\ref{fatt blocco H_{KL}}) and~(\ref{factrecur})
\begin{eqnarray}
\mathbf{H}^Q_{nn}            &=& \mathbf \Pi^Q_n \, \mathbf{\Gamma}^Q_n, \label{ideal1}\\
\mathbf{\Gamma}_{n}          &=& \mathbf{M}\mathbf{\Gamma}_{(n-1)},      \label{ideal2} 
\end{eqnarray}
This can be considered as an ideal algorithm. Feeding into the
system~(\ref{ideal1}), (\ref{ideal2}) the input
$\mathbf{H}^Q_{nn}$, which is known since
$Q$ is given, produces the output $\mathbf{M}$, whose blocks
contain the parameters $M(y)$ sought for. In real situations these
exact factorizations are generally not valid since $Q$ might not be a HMM, or might be
a HMM of order larger than $N$. This suggests constructing a two
step algorithm where~(\ref{ideal1}) and~(\ref{ideal2})
are substituted with NMFs. The scheme below
illustrates the two steps.

\medskip
{\it Step 1: Law Approximation Step}
\noindent\begin{eqnarray*}
\mbox{Given:} & \mathbf H^Q_{nn} \qquad\qquad\qquad\qquad\qquad\quad\\
\mbox{Problem:} & \underset{\mathbf \Pi_n,\mathbf{\Gamma}_n}\min D(\mathbf H^Q_{nn}\|\boldsymbol\Pi_n \mathbf{\Gamma}_n)\qquad\qquad\;\\
 & \mbox{constraints } e^\top \Pibf_n e=1, \boldsymbol{\Gamma}_n e=e\\
\mbox{Solution:} & \mathbf \Pi_{\,n}^*,\,\, \mathbf{\Gamma}_n^*.\qquad\qquad\qquad\qquad\quad
\end{eqnarray*}
Note that $\mathbf \Pi_{\,n}^*$ and $\mathbf{\Gamma}_n^*$ are of respective sizes $(m^n\times N)$ and $(N\times m^n)$.
Using equation~(\ref{newgamma}) one can build the matrix $\mathbf{\Gamma}_{(n-1)}^*$ needed below.

\bigskip
{\it Step 2: Parametrization Step (version with $\bf\Gamma$)}
\noindent\begin{eqnarray*}
\mbox{Given:} & \mathbf{\Gamma}_n^*, \,\, \mathbf{\Gamma}_{(n-1)}^* \qquad\qquad\qquad\qquad\quad\\
\mbox{Problem:} & \underset{\mathbf M}\min\, D\left(\mathbf{\Gamma}_{n}^{*}\|\mathbf{M\Gamma}_{(n-1)}^*\right)\qquad\qquad \\
  & \mbox{constraint } \mathbf Me=e \qquad\qquad\qquad \\
 \mbox{Solution:}& \mathbf{M}^{*}=\left[M^{*}(y_1)\ldots M^{*}(y_m)\right].\qquad
\end{eqnarray*}

\smallskip\noindent
Note that the constraint $\mathbf Me=e$, imposed at Step~2,
corresponds to the requirement that the transition matrix of the
underlying Markov chain be stochastic. The resulting
$A^*=\sum_{y}M^{*}(y)$ can be used to compute the parameter $\pi^*= \pi^*A^*$.

\smallskip Step~1 of the algorithm behaves like a typical EM method, with convergence of the
NMF algorithm to local minima of the divergence and dependence on the initial conditions $\mathbf \Pi_n^0,\mathbf{\Gamma}_n^0$.
Step 2 behaves much better, as one of the factors is fixed. The following Lemma summarizes the convergence properties of NMF
algorithms in the general setup of Step~2 of the algorithm. Although an easy consequence of the results in~\cite{ct1984}, to the best of our
knowledge the Lemma was not introduced before, at least not in the literature on NMF method
\smallskip
\begin{lemma} Let $S$ and $\Gamma$ be given stochastic matrices of sizes $m \times m$ and $N \times m$ respectively.
Let $M_0$ be a given stochastic matrix of size $m\times N$ with strictly positive elements. The NMF iterative algorithm~\cite{lee},~\cite{finessospreij2005} applied to the problem
\begin{equation} \label{csiszar-matrix}
\underset{M}\min D(S \|M \Gamma)
\end{equation}
with initial condition $M_0$, produces a sequence of matrices $M_k \rightarrow M^*$ (elementwise) where $M^*$ is the minimizer
of $D(S \|M \Gamma)$.
\end{lemma}
\smallskip\noindent The proof follows directly from Theorem 5 of~\cite{ct1984}, once it is recognized that the problem is decoupled in the rows of $M$.

\medskip
The analysis of the equations~(\ref{fatt blocco H_{KL}}) and~(\ref{factrecur}), valid in the ideal HMM case,
shows that as an alternative one could write the system
\begin{eqnarray*}
\mathbf{H}^Q_{nn}            &=& \mathbf \Pi^Q_n \, \mathbf{\Gamma}^Q_n, \\
\mathbf{\widetilde{\Pi}}_{n} &=& \mathbf{\Pi}_{n-1} \mathbf{M}.
\end{eqnarray*}

\noindent This suggests that it would be possible to substitute Step~2 of the proposed algorithm with the alternative Step2.alt below.
Note that from the output $\mathbf \Pi_{\,n}^*$ of Step~1, one can construct, as noted in Remark~\ref{remark1}, the matrix
$\mathbf \Pi_{\,n-1}^*$ and also repackage it, using equation~(\ref{factrecur}), to form the matrix $\mathbf{\tilde{\Pi}}_{n}^{*}$.

\medskip
{\it Step 2.alt: Parametrization Step (version with $\bf\Pi$)}
\noindent\begin{eqnarray*}
\mbox{Given:} & \mathbf{\tilde{\Pi}}_{n}^{*},\,\, \mathbf{\Pi}_{n-1}^{*} \qquad\qquad\qquad\qquad\quad\\
\mbox{Problem:} & \underset{\mathbf M}\min\, D\left(\mathbf{\tilde{\Pi}}_{n}^{*}\|\mathbf{\Pi}_{n-1}^{*} \mathbf{M}\right) \qquad\qquad\\
& \mbox{constraint } \mathbf Me=e  \qquad\qquad\qquad\\\
 \mbox{Solution:}& \mathbf{M}^{*}=\left[M^{*}(y_1)\ldots M^{*}(y_m)\right].\qquad
 \end{eqnarray*}

%
\section{Numerical Simulations}
In this section we show the results of numerical simulations designed to assess the performance
of the algorithm on HMM order reduction problems. We have also tried to evaluate possible
differences in performance related to the version of Step~2 being used.
We selected two examples. In both cases the given process is a 4 state HMM with binary output.

\smallskip
{\it Example 1:} The given HMM is specified by its transition probabilities matrix $A$ and
its read-out matrix $B$. Recall that $M(y)=AB_y$ where $B_y$ is the diagonal matrix with the $y$-th column of $B$ on the diagonal.
For this example $y\in\{0,1\}$.

\[A_1=\left(
      \begin{array}{cccc}
        0.3   & 0.15  & 0.1   & 0.45 \\
        0.1   & 0.5   & 0.2   & 0.2 \\
        0.25  & 0.15  & 0.35  & 0.25 \\
        0.2   & 0.35  & 0.4   & 0.05 \\
      \end{array}
    \right),
\]

\[B=\left(
      \begin{array}{cc}
        0.3   & 0.7 \\
        0.4   & 0.6 \\
        0.9   & 0.1 \\
        0.7   & 0.3 \\
      \end{array}
    \right).
\]

\smallskip
{\it Example 2:} Same read-out matrix $B$ as Example 1. The transition probabilties matrix is

\[A_2=\left(
      \begin{array}{cccc}
        0.125 & 0.3   & 0.025 & 0.55 \\
        0.4   & 0.4   & 0.025 & 0.175 \\
        0.\bar{3} & 0.\bar{3} & 0 & 0.\bar{3} \\
        0.2   & 0.5   & 0.025 & 0.275 \\
      \end{array}
    \right).
\]

\subsection{Discussion of Step 2 of the Algorithm}
\begin{figure}[bh]
\centering
    \includegraphics[bb=80.64 59.04 681.12 400.32,width=5.3cm, height=3.2cm]{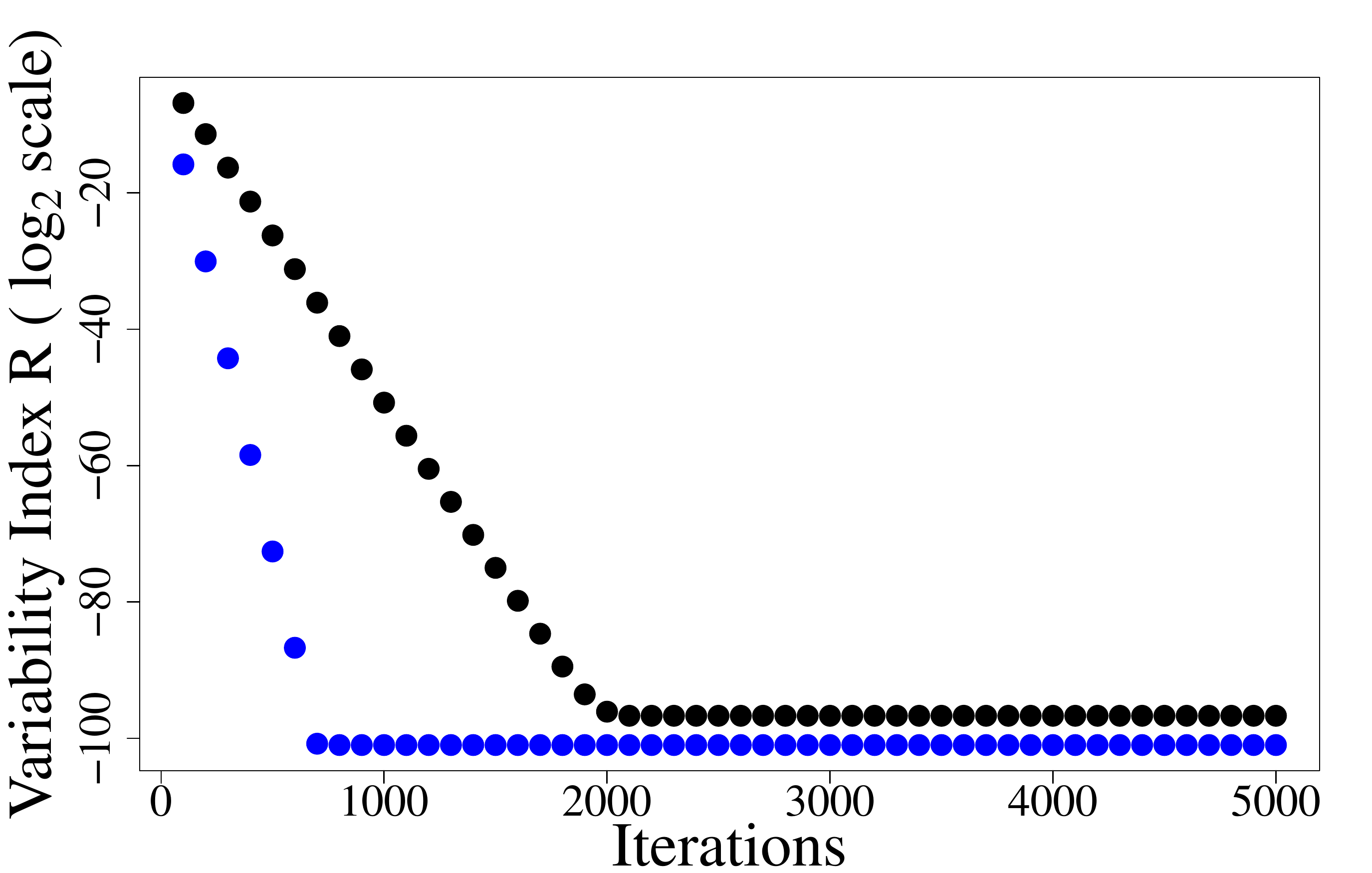} \hspace{.5em}\qquad\\\vspace{1em}
 \caption{Example 1 -- Variability of $\bf M^*$ in Step 2 (average of 30 runs), \mbox{$\bf\Gamma$ in black, $\bf\Pi$ in blue.}} \label{fig:ex3.variabilityM.N4NP2}
\end{figure}

\begin{figure}[h]
\centering
    \includegraphics[bb=80.64 59.04 637.92 400.32,width=5cm, height=3.2cm]{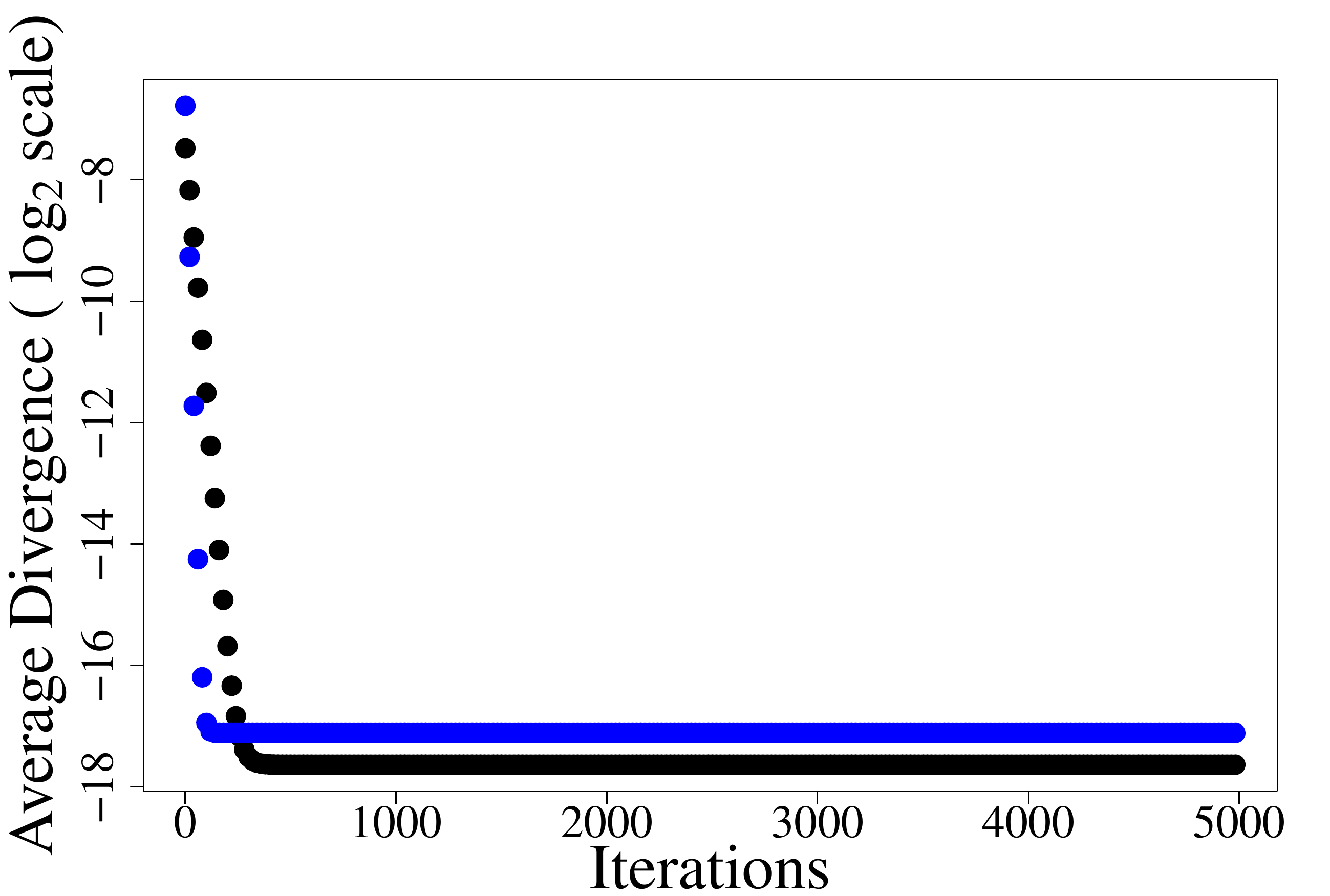} \hspace{.5em}\qquad\\\vspace{1em}
 \caption{Example 1 -- Divergence decay in Step 2 (average of 30 runs), \mbox{$\bf\Gamma$ in black, $\bf\Pi$ in blue.~~~~~~~~~~~~~~~~~~~~~~~~~~~~~~~~~~~~~~~~~~~~~~~~~~~~~~~~~~~~} \mbox{~~~~~~~~~~~~~~~~~~~~~~~~~~~~~~~~~~~~~~~~~~~~~~~~~~~~~}} \label{fig:ex3.I-divergence.N4NP2}
\end{figure}

\begin{figure}[ht]
\centering
    \includegraphics[bb=80.64 59.04 681.12 400.32,width=5.3cm, height=3.2cm]{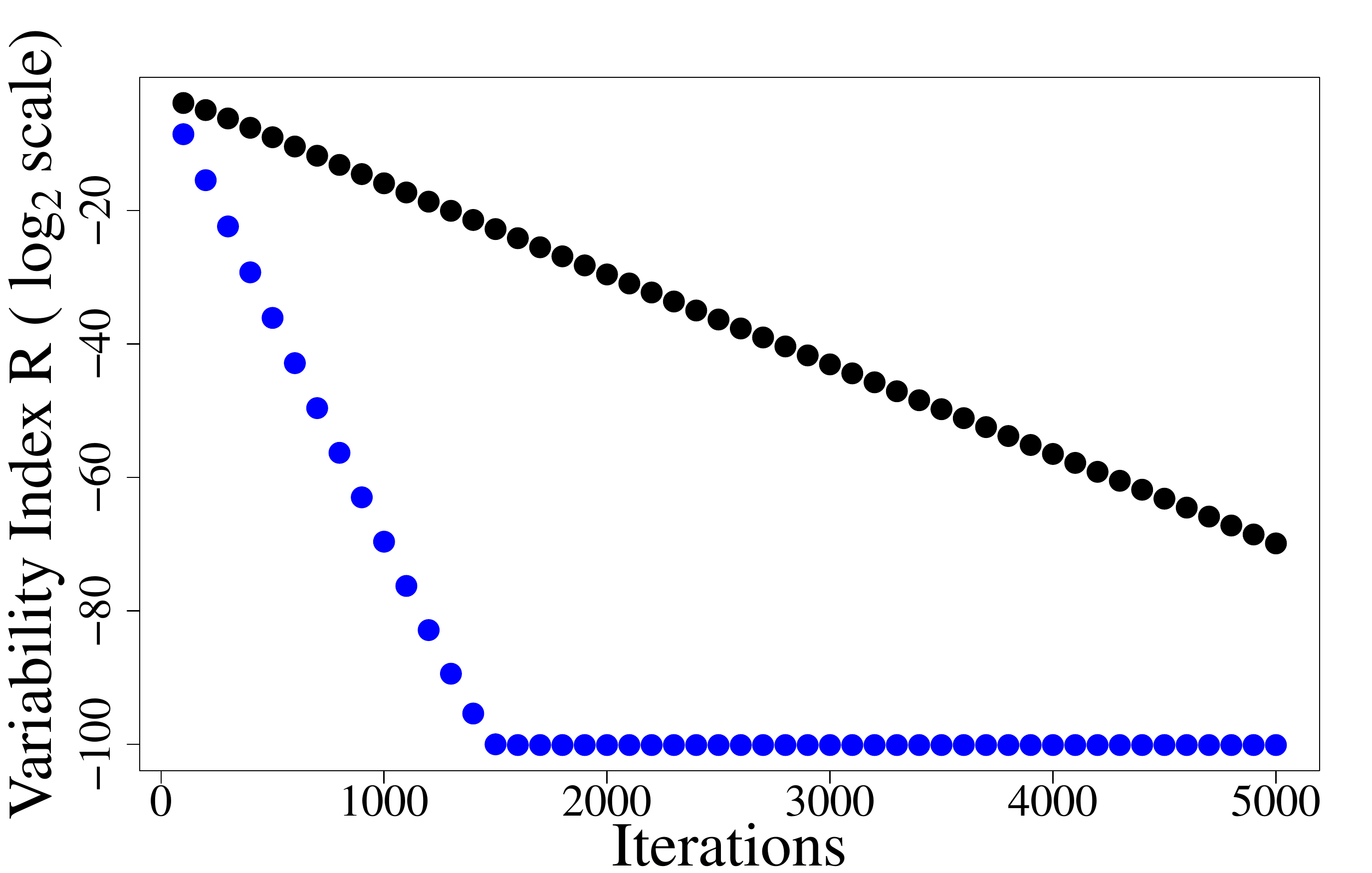} \hspace{.5em}\qquad\\\vspace{1em}
 \caption{Example 2 -- Variability of $\bf M^*$ in Step 2 (average of 30 runs), \mbox{$\bf\Gamma$ in black, $\bf\Pi$ in blue.}} \label{fig:ex4.variabilityM.N4NP2}
\end{figure}

\begin{figure}[h]
\centering
    \includegraphics[bb=80.64 59.04 637.92 400.32,width=5cm, height=3.2cm]{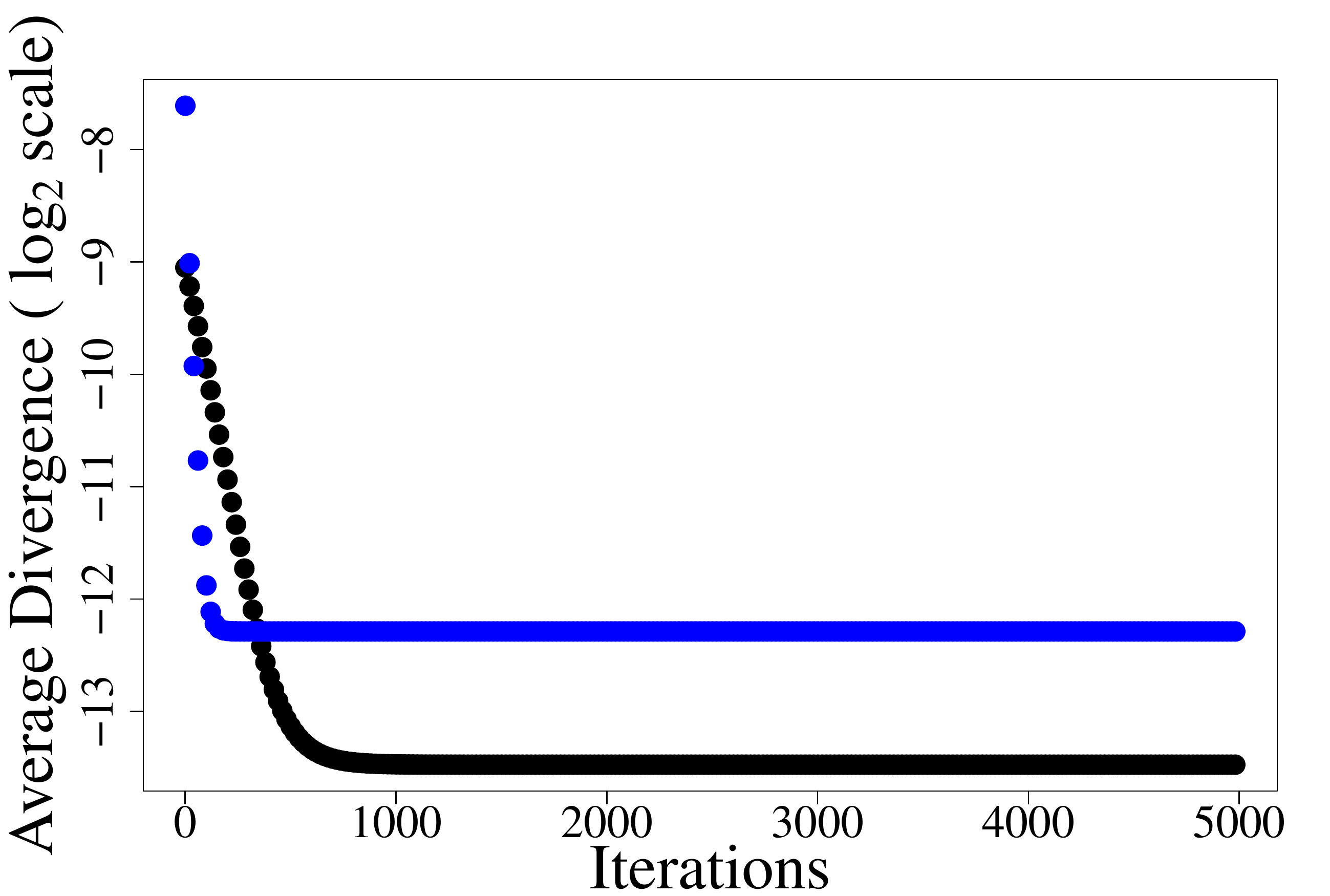} \hspace{.5em}\qquad\\\vspace{1em}
 \caption{Example 2 -- Divergence decay in Step 2 (average of 30 runs), \mbox{$\bf\Gamma$ in black, $\bf\Pi$ in blue.}} \label{fig:ex4.I-divergence.N4NP2}
\end{figure}

As seen at the end of the previous section, Step~2 of the algorithm
can be implemented in two versions. We will call the first the $\Gamma$
version and the alternative the $\Pi$ version.
By Lemma 1, and a variant of it, both NMFs are guaranteed to converge
to the global optimum, irrespective of the initial condition $M^0$.
We are interested in evaluating numerically if there are significant
differences in the speed of convergence of the $\Gamma$ and $\Pi$ versions
of Step~2. For the purpose of this comparison Step~1 was fixed and taken
as an order reduction step, from size 4 to size 2, for both Example~1 and Example~2.
We tested, on each example, the $\Gamma$ and the $\Pi$ versions of Step~2,
running for each 30 NMFs choosing $ M^0$ randomly.
The results were compared in terms of speed of convergence of the divergence
and variability of the resulting $\bf M^*$'s.
As an index of variability of the $\bf M^*$ parameters we computed the sum
of the sample variances of the elements of $\bf M^{*}$ resulting for each
of the 30 runs, i.e.
\begin{equation*}
    R=\sum_{i=1}^N\sum_{j=1}^{2N}\frac{1}{T-1}\sum_{t=1}^T(M_{ij}^{*t}-\bar{M}_{ij}^*)^2,
\end{equation*}
where $N=2$, is the size of the reduced HMM, $T=30$ is the number of different initial conditions $M^0$,
and  $\bar{M}_{ij}^*=\frac{1}{T}\sum_{t=1}^T M_{ij}^{*t}$.
Fig.~\ref{fig:ex3.variabilityM.N4NP2} and Fig.~\ref{fig:ex3.I-divergence.N4NP2}, relative to Example~1, display the variability index $R$ and the averaged (over the 30 runs with randomly chosen $M^0$) divergence decay against the number of NMF iterations. Fig.~\ref{fig:ex4.variabilityM.N4NP2} and Fig.~\ref{fig:ex4.I-divergence.N4NP2} show the same for Example~2.

As expected, both versions of Step~2 are not sensitive to the initial
conditions and the NMF iterations converge. What is particularly pleasing is that the resulting $\bf M^{*}$ is the same for both versions
of Step~2. The divergence appears to converge much faster than the parameter $\bf M^{*}$, a clear sign of the flatness of the criterion near optimality. Indicating with $\bf M^{*}_\Gamma$ and $\bf M^{*}_\Pi$ the parameters of the best approximation obtained with the two versions of Step 2,
we see that, for Example 1
\[\bf M^{*}_\Gamma=\left(
                    \begin{array}{cccc}
                      0.26633 & 0.24781 & 0.31323 & 0.17263 \\
                      0.20449 & 0.44350 & 0.21345 & 0.13856 \\
                    \end{array}
                  \right),
\]
\[\bf M^{*}_\Pi=\left(
                    \begin{array}{cccc}
                      0.26633 & 0.24778 & 0.31324 & 0.17265 \\
                      0.20453 & 0.44349 & 0.21344 & 0.13855 \\
                    \end{array}
                  \right),
\]
\mbox{~~~~~~~~}\\

\vspace{15em}
\vspace{1em}

and for Example 2
\[\bf M^{*}_\Gamma=\left(
                    \begin{array}{cccc}
                    0.28811 & 0.15352 & 0.29598 & 0.26239 \\
                    0.25952 & 0.25451 & 0.14339 & 0.34258 \\
                    \end{array}
                  \right),
\]
\[\bf M^{*}_\Pi=\left(
                    \begin{array}{cccc}
                     0.28814 & 0.15364 & 0.29571 & 0.26251 \\
                     0.25937 & 0.25451 & 0.14342 & 0.34271 \\
                    \end{array}
                  \right).
\]

Results along the same lines were obtained running many examples, corresponding to different sizes of order reduction, and/or different given HMM processes. The speeds of convergence of the divergence for the two versions of Step~2 remain always comparable. The variability of
$\bf M^{*}$ is sometimes slightly different for the two versions. The conclusion we drew is that there is not a version that systematically outperforms the other. For top performance in terms of variability it seems more appropriate to evaluate on a case by case basis which
version is most suitable.

%
%

\subsection{HMMs order reduction}
We collect here the results of four experiments of HMM order reduction. For each of the HMMs of size 4 given as Example~1 and Example~2, we ran the algorithm to produce HMMs of order reduced to 3 and 2. In all cases Step~2 was performed in the $\bf\Gamma$ version.
The set-up of the algorithm was fixed as follows: 3,000 NMF iterations for Step~1 in all cases, 3,000 and 20,000 NMF iterations for Step~2 for the order reductions $4 \rightarrow 2$ and $4 \rightarrow 3$ respectively. The number of NMF iterations was decided fixing a threshold on the $M^*$ variability index, $R$.
We randomly chose 30 initial conditions $(\mathbf \Pi_{\,n}^0,\,\, \mathbf{\Gamma}_n^0)$ for Step~1. As the NMF problem posed in Step~1 is sensitive to initial conditions (it behaves like a standard EM) we obtained 30 different pairs $(\bf\Pi_n^*,\Gamma_n^*)$ at its output.
We therefore had 30 different NMF problems as input to Step~2 leading to different $\bf M^*$.
In the following tables we report the divergence values, before and after each of the
two steps of the algorithm (DIV1b, DIV1, DIV2b, DIV2) and the final divergence (DIV) between the original process and the best approximation
obtained with the resulting $\bf M^*$. Specifically Table~\ref{tab:example3_N4_NP2} and Table~\ref{tab:example4_N4_NP2} refer to the HMM
order reduction $4 \rightarrow 2$ for Example~1 and Example~2 respectively. Likewise Table~\ref{tab:example3_N4_NP3} and Table~\ref{tab:example4_N4_NP3}
report the results of the order reduction $4 \rightarrow 3$ for the two examples.
The  high variability of the divergence after Step~2 (DIV2), especially in the case $4 \rightarrow 3$, is not surprising.
It depends on the fact that each run corresponds to a different NMF problem for Step~2.

\begin{figure}[ht]
\centering
    \includegraphics[bb=87.84 66.24 673.92 493.92,width=5.5cm, height=4.25cm]{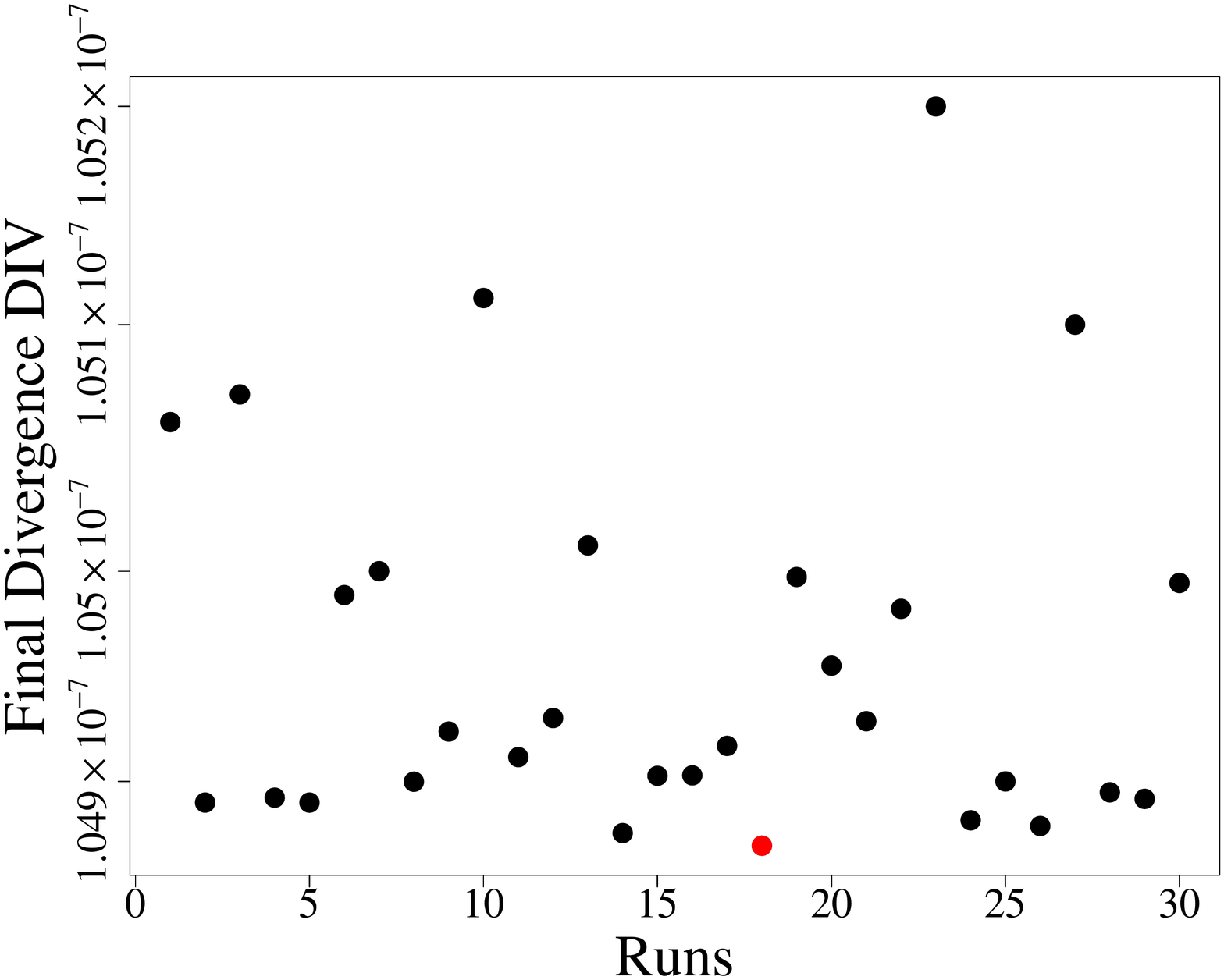} \hspace{.5em}\qquad\\\vspace{1.5em}
 \caption{Example 1 -- reduction $4\rightarrow 2$ -- final divergence of the 30 runs.} \label{fig:ex3.divergence.N4NP2}
\end{figure}

\begin{figure}[ht]
\centering
    \includegraphics[bb=87.84 66.24 673.92 493.92,width=5.5cm, height=4.25cm]{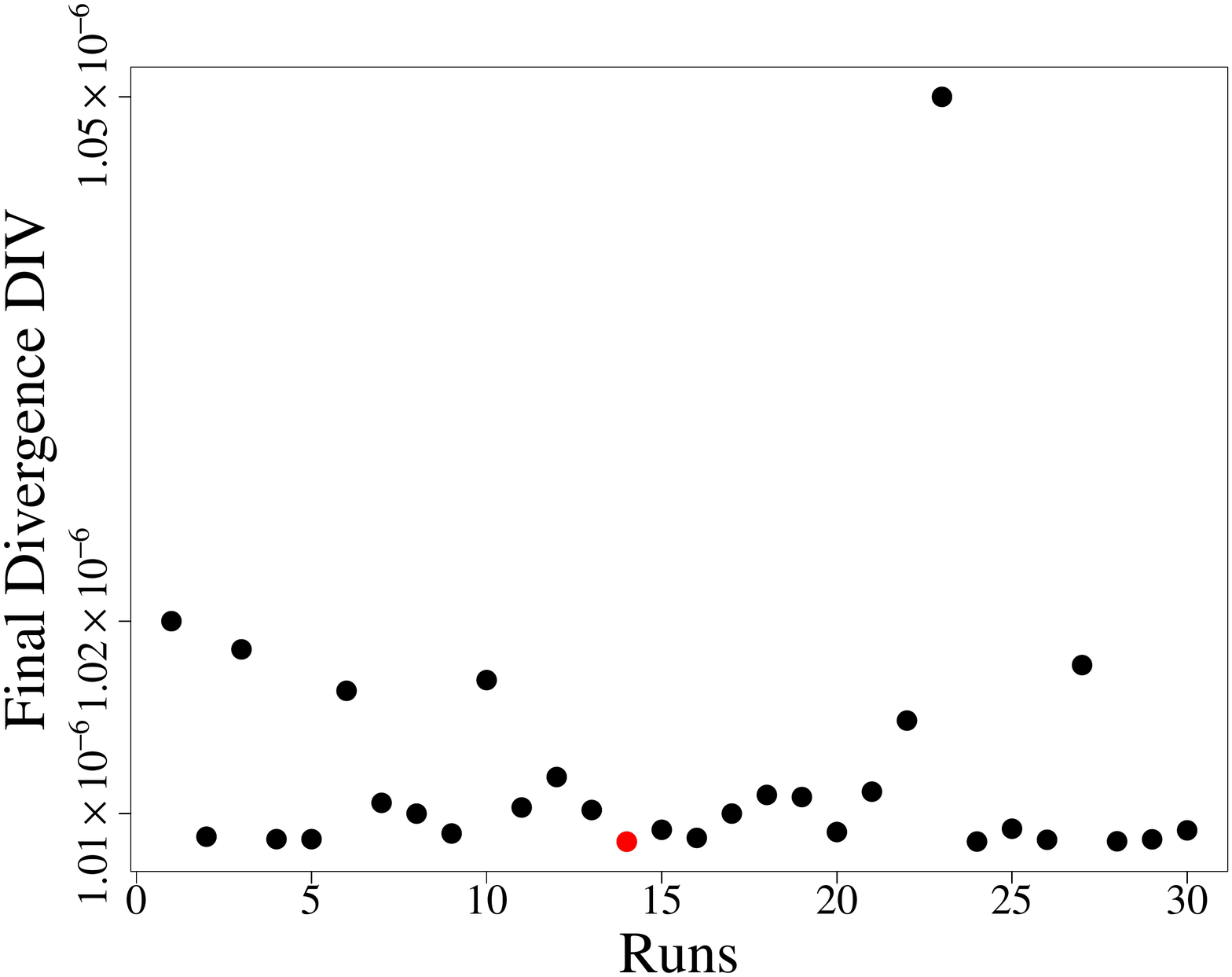} \hspace{.5em}\qquad\\\vspace{1.5em}
 \caption{Example 2 -- reduction $4\rightarrow 2$ -- final divergence of the 30 runs.} \label{fig:ex4.divergence.N4NP2}
\end{figure}

\begin{figure}[h]
\centering
    \includegraphics[bb=87.84 66.24 673.92 493.92,width=5.5cm, height=4.25cm]{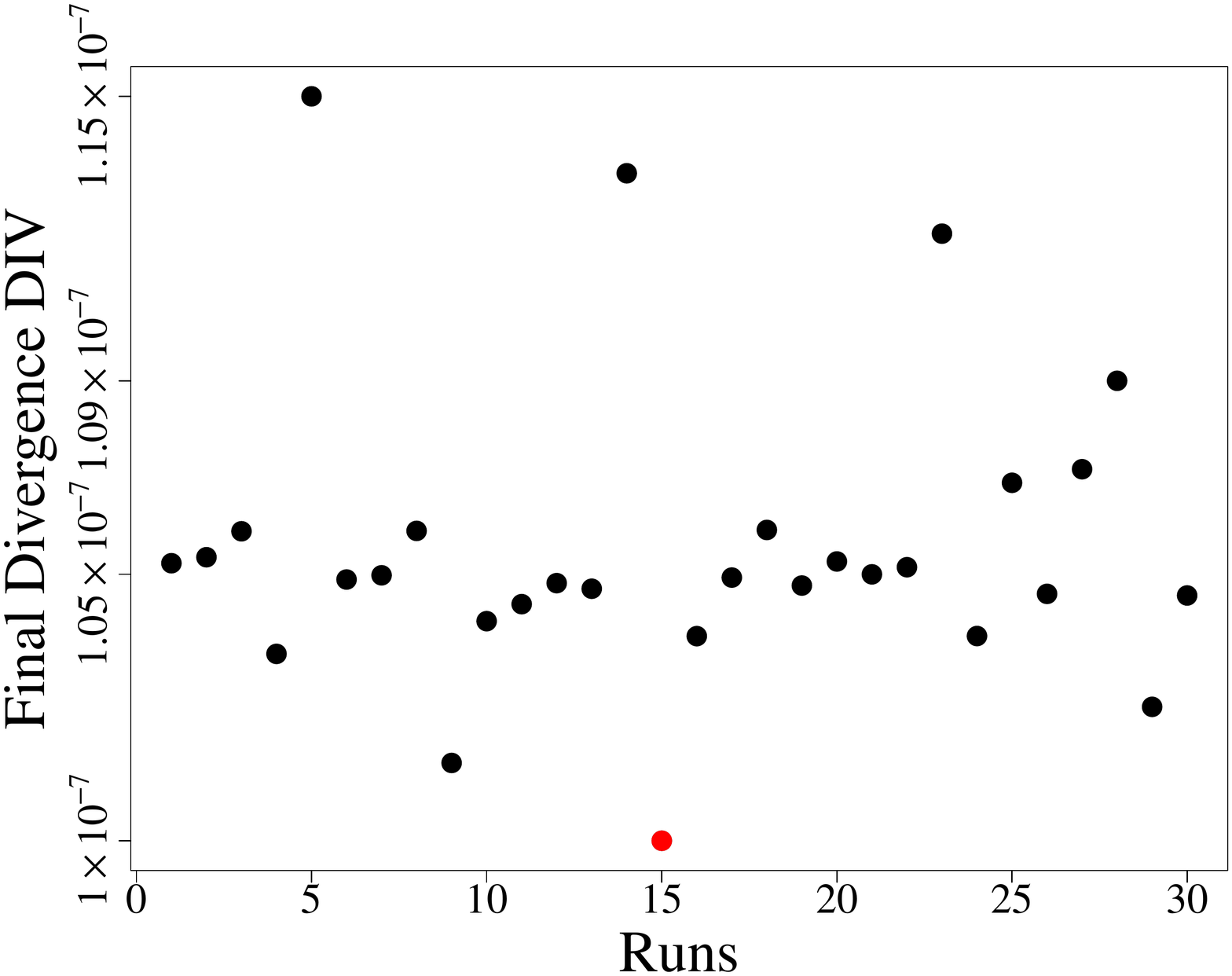} \hspace{.5em}\qquad\\\vspace{1.5em}
 \caption{Example 1 -- reduction $4\rightarrow 3$ -- final divergence of the 30 runs.} \label{fig:ex3.divergence.N4NP3}
\end{figure}

\begin{figure}[h]
\centering
    \includegraphics[bb=87.84 66.24 673.92 493.92,width=5.5cm, height=4.25cm]{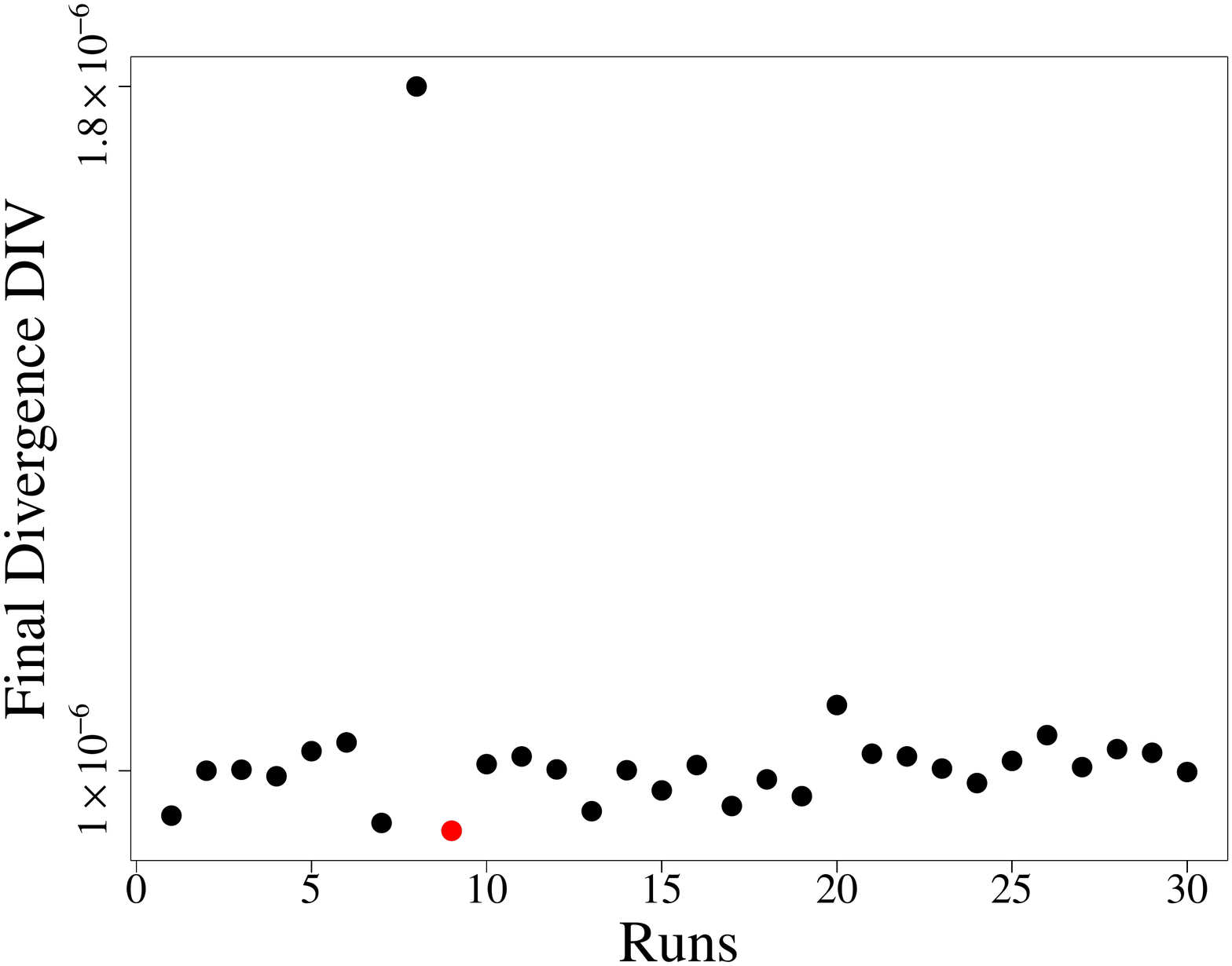} \hspace{.5em}\qquad\\\vspace{1.5em}
 \caption{Example 2 -- reduction $4\rightarrow 3$ -- final divergence of the 30 runs.} \label{fig:ex4.divergence.N4NP3}
\end{figure}

\smallskip It is interesting to notice that the final divergence (DIV) between the original process and the best approximation have the same order of magnitude, irrespective of the random initializations for the NMFs of the two steps.
In Fig.~\ref{fig:ex3.divergence.N4NP2}, relative to Example~1, and Fig.~\ref{fig:ex4.divergence.N4NP2}, relative to Example~2,
are plotted the divergences between the original process and the output process of the algorithm for the order reduction
case $4 \rightarrow 2$. Fig.~\ref{fig:ex3.divergence.N4NP3} and Fig.~\ref{fig:ex4.divergence.N4NP3} show the case $4 \rightarrow 3$.
In the plots the best numerical approximations are highlighted in red.\vspace{.1em}

%

\begin{table}[h]
  \centering
  \caption{Example 1  -- order reduction $4 \rightarrow 2$ -- divergence decays}
    \begin{tabular}{rrrrrr}
    \addlinespace
    \toprule
    RUN  & DIV1b & DIV1$\times 10^{9}$ & DIV2b & DIV2$\times 10^{6}$  & DIV$\times 10^{7}$ \\
    \midrule
    1     & 0.027 & 6.470 & 0.031 & 11.69 & 1.051 \\
    2     & 0.027 & 6.470 & 1.271 & 4.920 & 1.049 \\
    3     & 0.026 & 6.470 & 0.221 & 12.08 & 1.051 \\
    4     & 0.026 & 6.470 & 0.092 & 5.088 & 1.049 \\
    5     & 0.027 & 6.470 & 0.066 & 4.918 & 1.049 \\
    6     & 0.027 & 6.470 & 0.611 & 9.308 & 1.050 \\
    7     & 0.027 & 6.470 & 0.562 & 9.844 & 1.050 \\
    8     & 0.026 & 6.470 & 0.231 & 5.522 & 1.049 \\
    9     & 0.025 & 6.470 & 0.020 & 6.803 & 1.049 \\
    10    & 0.028 & 6.470 & 0.565 & 13.23 & 1.051 \\
    11    & 0.026 & 6.470 & 0.048 & 6.136 & 1.049 \\
    12    & 0.027 & 6.470 & 0.243 & 7.058 & 1.049 \\
    13    & 0.027 & 6.470 & 0.108 & 10.09 & 1.050 \\
    14    & 0.028 & 6.470 & 0.347 & 3.833 & 1.049 \\
    15    & 0.026 & 6.470 & 0.378 & 5.667 & 1.049 \\
    \bottomrule
    \end{tabular}
  \label{tab:example3_N4_NP2}
\end{table}
%
%
\begin{table}[htbp]
  \centering
  \caption{Example 2  -- order reduction $4 \rightarrow 2$ -- divergence decays}
    \begin{tabular}{rrrrrr}
    \addlinespace
    \toprule
    RUN  & DIV1b & DIV1$\times 10^{8}$  & DIV2b & DIV2$\times 10^{4}$  & DIV$\times 10^{6}$ \\
    \midrule
    1     & 0.015 & 6.601 & 0.046 & 4.198 & 1.020 \\
    2     & 0.014 & 6.601 & 0.886 & 0.878 & 1.007 \\
    3     & 0.014 & 6.601 & 0.154 & 3.915 & 1.018 \\
    4     & 0.014 & 6.601 & 0.023 & 0.744 & 1.007 \\
    5     & 0.014 & 6.601 & 0.009 & 0.737 & 1.007 \\
    6     & 0.013 & 6.601 & 0.212 & 3.468 & 1.016 \\
    7     & 0.013 & 6.601 & 0.452 & 1.814 & 1.009 \\
    8     & 0.013 & 6.601 & 0.137 & 1.612 & 1.008 \\
    9     & 0.014 & 6.601 & 0.076 & 1.014 & 1.007 \\
    10    & 0.014 & 6.601 & 0.287 & 3.579 & 1.017 \\
    11    & 0.013 & 6.601 & 0.154 & 1.761 & 1.009 \\
    12    & 0.014 & 6.601 & 0.219 & 2.330 & 1.011 \\
    13    & 0.014 & 6.601 & 0.364 & 1.666 & 1.009 \\
    14    & 0.015 & 6.606 & 0.273 & 0.496 & 1.007 \\
    15    & 0.014 & 6.601 & 0.215 & 1.148 & 1.007 \\
    \bottomrule
    \end{tabular}
  \label{tab:example4_N4_NP2}
\end{table}
%
%
\begin{table}[htbp]
  \centering
  \caption{Example 1  -- order reduction $4 \rightarrow 3$ -- divergence decays}
    \begin{tabular}{rrrrrr}
    \addlinespace
    \toprule
    RUN  & DIV1b & DIV1$\times 10^{9}$  & DIV2b & DIV2$\times 10^{3}$  & DIV$\times 10^{7}$ \\
    \midrule
    1     & 0.022 & 6.468 & 0.126 & 28.306 & 1.052 \\
    2     & 0.022 & 6.483 & 0.033 & 4.616 & 1.053 \\
    3     & 0.022 & 6.470 & 0.253 & 0.022 & 1.059 \\
    4     & 0.021 & 11.444 & 0.055 & 0.493 & 1.034 \\
    5     & 0.021 & 6.438 & 0.419 & 0.015 & 1.146 \\
    6     & 0.021 & 6.469 & 0.036 & 1.693 & 1.049 \\
    7     & 0.021 & 6.477 & 0.105 & 3.387 & 1.050 \\
    8     & 0.022 & 6.454 & 0.237 & 0.030 & 1.059 \\
    9     & 0.022 & 6.451 & 0.134 & 0.017 & 1.012 \\
    10    & 0.023 & 6.452 & 0.352 & 19.718 & 1.041 \\
    11    & 0.022 & 6.467 & 0.228 & 0.018 & 1.044 \\
    12    & 0.022 & 6.467 & 0.462 & 14.549 & 1.048 \\
    13    & 0.022 & 6.463 & 0.294 & 7.871 & 1.047 \\
    14    & 0.022 & 6.391 & 0.158 & 0.017 & 1.130 \\
    15    & 0.022 & 6.427 & 0.228 & 0.016 & 0.997 \\
    \bottomrule
    \end{tabular}
  \label{tab:example3_N4_NP3}
\end{table}
%
\begin{table}[htbp]
  \centering
  \caption{Example 2  -- order reduction $4 \rightarrow 3$ -- divergence decays}
    \begin{tabular}{rrrrrr}
    \addlinespace
    \toprule
    RUN  & DIV1b & DIV1$\times 10^{8}$  & DIV2b & DIV2$\times 10^{2}$  & DIV$\times 10^{6}$\\
    \midrule
    1     & 0.009 & 6.490 & 0.017 & 0.615 & 0.954 \\
    2     & 0.009 & 6.597 & 0.156 & 3.790 & 1.006 \\
    3     & 0.010 & 6.598 & 0.230 & 2.250 & 1.007 \\
    4     & 0.009 & 6.592 & 0.104 & 0.0316& 0.999 \\
    5     & 0.009 & 6.634 & 0.129 & 1.086 & 1.029 \\
    6     & 0.008 & 6.615 & 0.216 & 0.057 & 1.039 \\
    7     & 0.008 & 6.479 & 0.278 & 0.7781 & 0.945 \\
    8     & 0.009 & 7.738 & 0.120 & 2.458 & 1.803 \\
    9     & 0.009 & 6.576 & 0.165 & 0.073 & 0.936 \\
    10    & 0.009 & 6.598 & 0.069 & 0.783 & 1.013 \\
    11    & 0.008 & 6.607 & 0.099 & 3.472 & 1.022 \\
    12    & 0.009 & 6.601 & 0.198 & 5.010 & 1.007 \\
    13    & 0.009 & 6.568 & 0.122 & 0.052 & 0.959 \\
    14    & 0.009 & 6.601 & 0.213 & 4.282 & 1.006 \\
    15    & 0.009 & 6.552 & 0.114 & 2.329 & 0.983 \\
    \bottomrule
    \end{tabular}
  \label{tab:example4_N4_NP3}
\end{table}

\vspace{2.5em}

\vspace{6em}

\vspace{6em}

\newpage

%
%
%
%
%
%
%
%
%

\end{document}